\documentclass[12pt]{article}

\usepackage{amsmath}
\usepackage{wasysym}
\usepackage{latexsym}
\usepackage{amsfonts}
\usepackage{mathrsfs}
\usepackage{amssymb}
\usepackage{ifsym}
\usepackage{dsfont}
\usepackage[all]{xy}
\usepackage{authblk}
\usepackage{lipsum}

\newcommand\blfootnote[1]{%
  \begingroup
  \renewcommand\thefootnote{}\footnote{#1}%
  \addtocounter{footnote}{-1}%
  \endgroup
}

\usepackage{amsfonts}
\usepackage{amssymb}

\DeclareMathOperator{\Def}{Def}
\DeclareMathOperator{\sat}{Sat}

\def\D11{\Delta_1^1$-${\sf CA}_0}
\def\S11{\Sigma_1^1$-${\sf AC}_0}

\begin{document}

\title{The Barwise-Schlipf Theorem}
\author[1]{Ali Enayat}
\affil[1]{University of Gothenburg, Gothenburg, Sweden\newline
\texttt{ali.enayat@gu.se}}
\author[2]{James H. Schmerl}
\affil[2]{University of Connecticut, Storrs, CT 06269, USA\newline
\texttt{james.schmerl@uconn.edu}}
\maketitle

\begin{abstract}
In 1975 Barwise and Schlipf published a landmark paper whose main theorem
asserts that a nonstandard model $\mathcal{M}$ of $\mathsf{PA}$ (Peano arithmetic) is
recursively saturated iff $\mathcal{M}$ has an expansion that satisfies the
subsystem $\Delta _{1}^{1}$-$\mathsf{CA}_{0}$ of second order arithmetic. In
this paper we identify a crucial error in the Barwise-Schlipf proof of the
right-to-left direction of the theorem, and additionally, we offer a correct
proof of the problematic direction.
\end{abstract}

\blfootnote{We are grateful to Roman Kossak and Mateusz  Łełyk for their help in improving our exposition.}

In their seminal paper \cite{bs}, Barwise
and Schlipf initiated the study of recursively saturated models of \textsf{PA%
}\ with the following theorem.

\bigskip

\textsc{Theorem}: (Barwise-Schlipf \cite{bs}) \emph{If ${\mathcal{M}}
\models \mathsf{PA}$ is nonstandard, then the following are equivalent$:$ }

\smallskip

\emph{$(1)$ ${\mathcal{M}}$ is recursively saturated. }

\smallskip

\emph{$(2)$ There is ${\mathfrak{X}}$ such that $({\mathcal{M}},{\mathfrak{X}%
}) \models \D11$. }

\smallskip

\emph{$(3)$ $({\mathcal{M}}, \Def({\mathcal{M}})) \models \D11 + \S 11$.}

\bigskip

\noindent Their proof of $(1)\Longrightarrow (3)$ (\cite[Theorem~2.2]{bs})
uses Admissible Set Theory. In a reprise of this theorem by Smory\'{n}ski
\cite[Sect.~4]{smo81}, a more direct proof of this implication, attributed
to Feferman and Stavi (independently), is presented. This same proof is
essentially repeated by Simpson \cite[Lemma~IX.4.3]{sim}. The implication $%
(3)\Longrightarrow (2)$ is trivial. In the proof of the remaining
implication $(2)\Longrightarrow (1)$, it is claimed \mbox{%
\cite[Theorem~3.1]{bs}} that if ${\mathcal{M}}$ is nonstandard and not
recursively saturated and $\Def({\mathcal{M}})\subseteq {\mathfrak{X}}%
\subseteq {\mathcal{P}}(M)$, then $({\mathcal{M}},\mathfrak{X})\not\models \D%
11$ because the standard cut $\omega $ is $\Delta _{1}^{1}$-definable{%
\footnote{%
All usages of \emph{definable} in this paper should be understood as \emph{%
definable with parameters}.} in $({\mathcal{M}},{\mathfrak{X}})$. }To prove
that, they let $\Phi =\{\varphi _{n}(x):n<\omega \}$ be a finitely
realizable type that is not realized in ${\mathcal{M}}$, where $\langle
\varphi _{n}(x):n<\omega \rangle $ is a recursive sequence (with the
understanding that there is a finite set $F\subseteq M$ such that any
parameter occurring in any $\varphi _{n}(x)$ is in $F$). Then they let $%
Y=\{a_{m}:m<\omega \}$, where $a_{m}$ is the least $n<\omega $ such that ${%
\mathcal{M}}\models \lnot \varphi _{n}(m).$ Their $\Sigma _{1}^{1}$%
-definition of $Y$ is correct, but their purported $\Pi _{1}^{1}$-definition
of $Y$ does not work since it defines the set $Y\cup (M\setminus \omega )$.
Murawski's exposition of the Barwise-Schlipf theorem {\cite{m76} suffers
from the same gap. Smory\'{n}ski makes a similar error in his explicit claim
\cite[Lemma~4.2]{smo81} that if ${\mathcal{M}}$ is nonstandard and not
recursively saturated and $({\mathcal{M}},{\mathfrak{X}})\models \mathsf{ACA}%
_{0}$, then $\omega $ is $\Delta _{1}^{1}$-definable in $({\mathcal{M}},{%
\mathfrak{X}})$. We will show in Theorem 2 that this approach is doomed
since there are nonstandard models ${\mathcal{M}}$ that are not recursively
saturated even though $\omega $ is not $\Delta _{1}^{1}$-definable in $({%
\mathcal{M}},\Def({\mathcal{M}}))$. Nevertheless, we are still able to give
a proof (see Theorem~3) of $(2)\Longrightarrow (1)$. }

\smallskip

Suppose that ${\mathcal{M}}\models \mathsf{PA}$ and $A\subseteq M$. Then, $A$
is \textbf{recursively $\sigma $-definable} if there is a recursive sequence
$\langle \varphi _{n}(x):n<\omega \rangle $ of formulas (where for some
finite set $F\subseteq M$ any parameter occurring in any $\varphi _{n}(x)$
is in $F)$ such that each $\varphi _{n}(x)$ defines a subset $%
A_{n}\subseteq M$, with $A=\bigcup_{n<\omega }A_{n}$. For example, the
standard cut $\omega $ is recursively $\sigma $-definable, and so is any
finitely generated submodel of $\mathcal{M}$.

\bigskip

\textsc{Lemma 1}: \emph{Suppose that ${\mathcal{M}} \models \mathsf{PA}$ and
$A \subseteq M$. }

\smallskip

\emph{$(a)$ If $A$ is $\Sigma_1^1$-definable in $({\mathcal{M}}, \Def({%
\mathcal{M}}))$, then $A$ is recursively $\sigma$-definable. }

\smallskip

\emph{$(b)$ If ${\mathcal{M}}$ is not recursively saturated, $\Def({\mathcal{%
M}}) \subseteq {\mathfrak{X}} \subseteq {\mathcal{P}}(M)$ and $A$ is}

\hspace{14pt} \emph{recursively $\sigma $-definable, then $A$ is $\Sigma
_{1}^{1}$-definable in $({\mathcal{M}},{\mathfrak{X}})$.}\bigskip

\textit{Proof}.~(a) Suppose that $A$ is $\Sigma _{1}^{1}$-definable in $({%
\mathcal{M}},\Def({\mathcal{M}}))$ by the formula $\exists X\theta (x,X)$.
Let $\langle \psi _{n}(x,y):n<\omega \rangle $ be an enumeration of all
formulas in exactly two free variables, and let $\varphi _{n}(x,y)=\theta
(x,\{u:\psi_{n} (u,y)\}),$ i.e., $\varphi _{n}(x,y)$ is the result of substituting every occurrence of
subformulas of $\theta (x,X)$ of the form $t\in X$ (where $t$ is a term)
with $\varphi (t,y)$ (and re-naming variables to avoid unintended clashes).
Then $\langle \exists y\varphi _{n}(x,y):n<\omega \rangle $ is recursive,
each formula $\exists y\varphi _{n}(x,y)$ defines a subset $A_{n}$ of $A$, and $%
A=\bigcup_{n<\omega }A_{n}.$ Hence $A$ is recursively $\sigma $-definable.

\smallskip

(b)\hspace{1pt} Let $\sat(x,X)$ be a formula asserting that $X$ is a
satisfaction class for all formulas of length at most $x$. Let $A$ be
recursively $\sigma $-definable by the recursive sequence $\langle \varphi
_{n}(x):n<\omega \rangle $. We can assume that $\ell (\varphi _{n}(x))<\ell
(\varphi _{n+1}(x))$ for all $n<\omega $, where $\ell (\varphi (x))$ is the
length of $\varphi (x)$ (by replacing $\varphi _{n}(x)$ with $\bigvee_{i\leq
n}\varphi _{i}(x)$). The sequence $\langle \varphi _{n}(x):n<\omega \rangle $
is coded in ${\mathcal{M}}$, so let $d\in M$ be nonstandard such that $%
\langle \varphi _{n}(x):n<d\rangle $ extends $\langle \varphi
_{n}(x):n<\omega \rangle $ and $\ell (\varphi _{n}(x))$ is standard iff $n$
is. Then $A$ is $\Sigma _{1}^{1}$-definable in $({\mathcal{M}},\Def({%
\mathcal{M}}))$ by the formula $\exists X\theta (x,X)$, where
\begin{equation*}
\theta (x,X)=\exists z[\sat(z,X)\wedge \exists n<d\big(\ell (\varphi
_{n})\leq z\wedge \langle \varphi _{n},x\rangle \in X\big)].
\end{equation*}%
Thus, $A$ is $\Sigma _{1}^{1}$-definable in $({\mathcal{M}},\Def({\mathcal{M}%
}))$. The same definition works in $({\mathcal{M}},{\mathfrak{X}})$. \hfill $%
\square $

\bigskip

According to a definition in \cite[Notation 2.1(b)]{ks2}, the standard cut $%
\omega $ is said to be \textit{recursively definable} in a model $\mathcal{M}
$\ of $\mathsf{PA}$ if there is a recursive type $\Sigma (x)$
that is finitely realizable in $\mathcal{M}$, and which has the property
that for every elementary extension $\mathcal{N}$ of $\mathcal{M}$ that has
an element $b$ realizing $\Sigma (x)$, $b$ fills the standard
cut of $\mathcal{M}$, i.e., $b$ is a nonstandard element of $\mathcal{N}$
that is below all nonstandard elements of $\mathcal{M}$. \bigskip

\textsc{Lemma 2}: \emph{If }$\omega $\emph{\ is not recursively definable in
}$M$\emph{, then }$\omega $\emph{\ is not }$\Pi _{1}^{1}$\emph{-definable in
}$(\mathcal{M},\Def(\mathcal{M})).$\bigskip

\textit{Proof}.~If $\omega $ is not recursively definable in $\mathcal{M}$%
, then $M\backslash \omega $ is not recursively $\sigma $-definable. Thus by
Lemma~1(a), $\omega $ is not $\Pi _{1}^{1}$-definable in $({\mathcal{M}},\Def%
({\mathcal{M}}))$. \hfill $\square $

\bigskip

\textsc{Theorem 3}: \emph{Every consistent completion $T$ of $\mathsf{PA}$
has a countable nonstandard model }${\mathcal{M}}$\emph{\ that is not
recursively saturated and such that $\omega $ is not $\Delta _{1}^{1}$%
-definable in $({\mathcal{M}},\Def({\mathcal{M}})).$ Moreover, }$\mathcal{M}$%
\emph{\ can be arranged to be finitely generated.}

\bigskip

\textit{Proof}. By Lemma 2, the proof of the first assertion of Theorem 2 is
complete once we exhibit a nonstandard model $\mathcal{M}$ of \textsf{PA}
that is not recursively saturated, and in which $\omega $ is not recursively
definable. Recall that a model $\mathcal{M}$ of \textsf{PA} is \textit{short
recursively saturated} if $\mathcal{M}$ realizes every recursive \textit{%
short} type over $\mathcal{M}$, i.e., recursive types that include a formula
of the form $x<m$ for some $m\in M$. It can be readily verified that $\omega
$ is not recursively definable in any short recursively saturated model. On
the other hand, given a completion $T$\ of $\mathsf{PA}$, it is routine to
construct a nonstandard model $\mathcal{M}$ of $T$ that is short recursively
saturated but not saturated since if $\mathcal{N}$  is a recursively
saturated model of $T$, and $\mathcal{N}_{0}\prec \mathcal{N}$ is
nonstandard and finitely generated, then we can choose $\mathcal{M}$ to be
the submodel of $\mathcal{N}$ consisting all those $b\in N$ that are less
than some $a\in N_{0}.$ This makes it evident that $T$ has a countable
nonstandard model that is not recursively saturated and such that $\omega $
is not $\Delta _{1}^{1}$-definable in $(\mathcal{M},\Def(\mathcal{M}))$\emph{%
. }Finally, in order to establish the moreover clause, we note that
according to \cite[Corollary 2.8]{ks2} every consistent completion $T$ of
\textsf{PA} has a finitely generated ${\mathcal{M}}\models T$ such that $%
\omega $ is not recursively definable in $\mathcal{M}$. So by Lemma 2 we are
done. \hfill $\square $

\bigskip

\textsc{Theorem 4}: \emph{If ${\mathcal{M}}$ is nonstandard and $({\mathcal{M%
}},{\mathfrak{X}})\models \D11$, then ${\mathcal{M}}$ is recursively
saturated.}

\bigskip

\textit{Proof}. We will show that if ${\mathcal{M}}$ is nonstandard and not
recursively saturated and ${\mathfrak{X}} \subseteq {\mathcal{P}}(M)$, then $%
({\mathcal{M}}, {\mathfrak{X}}) \not\models \D11$. We can assume that $({%
\mathcal{M}}, {\mathfrak{X}}) \models \mathsf{ACA}_0$. There are two cases
depending on whether ${\mathcal{M}}$ is short or tall.

\smallskip

\emph{${\mathcal{M}}$ is short}: Let $a\in M$ be such that the elementary
submodel of ${\mathcal{M}}$ generated by $a$ is cofinal in ${\mathcal{M}}$.
Let $\langle \varphi _{n}(x):n<\omega \rangle $ be a recursive sequence of
formulas (with $a$ as the only parameter) such that $\varphi _{n}(x)$
defines $d_{n}\in M$, where $d_{n}$ is the least element $x$ above all
elements that are definable from $a$ by a formula of length at most $n$. Let
$D=\{d_{n}:n<\omega \}$. $D$ is unbounded in $\mathcal{M}$ and $d_{n}\leq d_{n+1}$ for all $%
n<\omega .$ Since $({\mathcal{M}},{\mathfrak{X}})\models \mathsf{ACA}_{0}$,
then $D\not\in {\mathfrak{X}}$ as otherwise $\omega \in {\mathfrak{X}}$.
Clearly, $D$ is recursively $\sigma $-definable; its complement also is
(using the recursive sequence $\langle \psi _{n}(x):n<\omega \rangle $,
where $\psi _{0}(x)$ is $x<d_{0}$ and $\psi _{n+1}(x)$ is $d_{n}<x<d_{n+1}$%
). By Lemma~1(b), $D$ is $\Delta _{1}^{1}$-definable in $({\mathcal{M}},{%
\mathfrak{X}})$.

\smallskip

\emph{${\mathcal{M}}$ is tall}: Since ${\mathcal{M}}$ is tall and not
recursively saturated, there is a recursive sequence $\langle \varphi
_{n}(x):n<\omega \rangle $ of formulas, among which is a formula $x<b$, that
is finitely realizable in ${\mathcal{M}}$ but not realizable in ${\mathcal{M}%
}$. According to \cite[Lemma~2.4]{ks1}\footnote{%
More explicitly, Lemma 2.4 of \cite{ks1} directly implies that if $\Sigma (x,%
\overline{a})$ is a recursive finitely realizable type over a model $%
\mathcal{M}$, and $\Sigma (x,\overline{a})$ includes $\{x<a_{0}\},$ then
there is a finitely realizable recursive type $\Gamma (x,\overline{a},d)$,
where $d$ can be chosen to be any nonstandard element of $\mathcal{M}$, such
that: (1) If $\mathcal{M}$ does not realize $\Sigma (x,\overline{a}),$ then $%
\mathcal{M}$ does not realize $\Gamma (x,\overline{a},d);$ (2) $\Gamma (x,%
\overline{a},d)$ consists of formula of the form $s_{n}(\overline{a},d)\leq x\leq
t_{n}(\overline{a},d)$, where both sequences $\left\langle s_{n}(\overline{x}%
,y):n<\omega \right\rangle $ and $\left\langle t_{n}(\overline{x}%
,y):n<\omega \right\rangle $ are recursive, and for all $n<\omega ,\
\mathcal{M}\models s_{n}(\overline{a},d)\leq s_{n+1}(\overline{a},d)\leq t_{n+1}(\overline{a},d)\leq
t_{n}(\overline{a},d).$}, we can assume that each $\varphi _{n}(x)$ defines an interval
$[a_{n},b_{n}]$, where $a_{n}\leq a_{n+1}\leq b_{n+1}\leq b_{n}$. Then, the
cut $I=\sup \{a_{n}:n<\omega \}=\inf \{b_{n}:n<\omega \}$, so both $I$ and
its complement are recursively $\sigma $-definable. Lemma~1 implies $I$ is $%
\Delta _{1}^{1}$-definable in $({\mathcal{M}},{\mathfrak{X}})$. Since $%
I\not\in {\mathfrak{X}}$, then $({\mathcal{M}},{\mathfrak{X}})\not\models \D%
11$. \hfill $\square $

\bigskip

We conclude with several remarks concerning the Theorem.

\smallskip

It is well known that $\Sigma _{k}^{1}$-$\mathsf{AC}_{0}$ implies $\Delta
_{k}^{1}$-$\mathsf{CA}_{0}$ for all $k<\omega $. An easy proof can be found
in \cite[Lemma~VII.6.6(1)]{sim}.

\smallskip

Barwise and Schlipf point out \cite[Remark, p.~52]{bs} that their
(erroneous) proof of $(2)\Longrightarrow (1)$ shows the slightly stronger
implication in which $\Delta _{1}^{1}$-$\mathsf{CA}_{0}$ is replaced by its
counterpart $\Delta _{1}^{1}$-$\mathsf{CA}_{0}^{-}$ in which there are no
set parameters. The same is true of our proof of $(2)\Longrightarrow (1)$.

\smallskip

The impression one might get from reading \cite{bs} is that $(1)\Longrightarrow
(3)$ is the deep direction of the Theorem (since it relies on the technology
of admissible sets), and  $(2)\Longrightarrow (1)$ is the fairly
straightforward one. Now, prospering from a 45 year hindsight, one can say that the exact opposite is the case: the hard direction is $%
(2)\Longrightarrow (1)$ while $(1)\Longrightarrow (3)$ can be handled by a
short proof based on first principles.



\begin{thebibliography}{9}
\bibitem{bs} Jon Barwise and John Schlipf, On recursively saturated models
of arithmetic, in: \emph{Model theory and algebra} (\textit{A memorial
tribute to Abraham Robinson}), Lecture Notes in Math., \textbf{498}, pp.
42--55, Springer, Berlin, 1975.

\bibitem{ks1} Matt Kaufmann and James H.~Schmerl, Saturation and simple
extensions of models of Peano arithmetic, Ann. Pure Appl. Logic \textbf{27}
(1984), 109--136.

\bibitem{ks2} Matt Kaufmann and James H.~Schmerl, Remarks on weak notions of
saturation in models of Peano arithmetic, J. Symbolic Logic \textbf{52}
(1987), 129--148.

\bibitem{m76} Roman Murawski, On expandability of models of Peano arithmetic
I, Studia Logica, \textbf{35} (1976) 409--419.

\bibitem{sim} Stephen G.~Simpson, \emph{Subsystems of second order arithmetic%
}, Perspectives in Mathematical Logic, Springer-Verlag, Berlin, 1999.

\bibitem{smo81} Craig Smory\'{n}ski, Recursively saturated nonstandard models
of arithmetic, J. Symbolic Logic \textbf{46} (1981), 259--286.
\end{thebibliography}
\end{document}